\documentstyle[12pt, righttag]
{amsart}

    \newtheorem{propo}{Proposition}[section]
    \newtheorem{theo}[propo]{Theorem}
    \newtheorem{lemma}[propo]{Lemma}

    \theoremstyle{definition}
    \newtheorem{defi}{Definition}[section]

    \theoremstyle{remark}

    \numberwithin{equation}{section}

	\newcommand{\nno}{\nonumber}
	\newcommand{\lbar}{\bigg\vert}
	
	\newcommand{\dps}{\displaystyle}

\newcommand{\wt}{\mbox{\rm wt}\ }
\renewcommand{\hom}{\mbox{\rm Hom}}

\begin{document}
\title[Introduction to Vertex operator algebras, III]{Introduction to
Vertex operator algebras, III}
\author[Yi-Zhi Huang]{Yi-Zhi Huang*}
\thanks{* Supported in part by NSF grant DMS-9301020 and
by DIMACS, an
NSF Science and Technology Center funded under contract STC-88-09648.}
\address{Department of Mathematics\\ Rutgers University\\
New Brunswick, NJ 08903}
\email{yzhuang@@math.rutgers.edu}

    \bibliographystyle{alpha}
    \maketitle

    \input amssym.def
    \input amssym

\tableofcontents

In this exposition, we continue the discussions of Dong \cite{D2}
and Li \cite{L}.
We shall prove an $S_{3}$-symmetry of the Jacobi identity, construct
the contragredient module for a module for a vertex operator algebra and apply
these to the construction of the vertex operator map for the moonshine module.
We shall introduce the notions of intertwining operator, fusion rule and
Verlinde algebra. We shall also describe briefly the geometric interpretation
of vertex operator algebras. We end the exposition with an explanation of
the role of vertex operator algebras in conformal field theories.

I would like to thank Masahiko Miyamoto for inviting me to this successful
conference and James Lepowsky for helpful mathematical comments.

\vspace{1em}
\noindent {\bf Notations}:

\noindent $\Bbb{C}$: the (structured set of) complex numbers.

\noindent $\Bbb{C}^{\times}$: the nonzero complex numbers.

\noindent $\Bbb{R}$: the real numbers.

\noindent $\Bbb{Z}$: the integers.

\noindent $\Bbb{Z}_{+}$: the positive integers.

\noindent $\Bbb{N}$: the nonnegative integers.

\section{$S_{3}$-symmetry of the Jacobi identity
and contragredient modules}

The results and constructions discussed in this section are all
natural {}from the axiomatic viewpoint. But they also have practical
uses in some very concrete problems. Before going into the detailed
discussions, let us first recall one of those problems.

One of the most important examples of vertex operator algebra is the
moonshine module constructed by Frenkel, Lepowsky and Meurman
\cite{FLM1} \cite{FLM2}. (See the introduction of \cite{FLM2} for
a historical discussion, including the important role of Borcherds'
announcement \cite{B}.) The construction can be briefly described as
follows: {}From the Leech lattice $\Lambda$, one can construct an
untwisted vertex operator algebra $V_{\Lambda}$.
 The automorphism
$\theta: \Lambda\to \Lambda$ defined by $\theta(x)=-x$ for any $x\in
\Lambda$ induces an automorphism of $V_{\Lambda}$
which is still denoted $\theta$. One can construct a unique irreducible
$\theta$-twisted module $V_{\Lambda}^{T}$ for $V_{\Lambda}$.
 The automorphism
$\theta: \Lambda\to \Lambda$ also induces an automorphism of $V_{\Lambda}$
and is also denoted $\theta$.
Let $V_{\Lambda}^{+}$ and
$(V_{\Lambda}^{T})^{+}$ be spaces of fixed points of $\theta$ in
$V_{\Lambda}$ and $V_{\Lambda}^{T}$, respectively. Then the moonshine
module is $V^{\natural}=V_{\Lambda}^{+}\oplus (V_{\Lambda}^{T})^{+}$
as a $\Bbb{Z}$-graded vector space. In \cite{FLM2}, the vertex
operator map for the moonshine module is defined and it is shown that
$V^{\natural}$ is indeed a vertex operator algebra.

The definition of vertex operator map for $V^{\natural}$ in
\cite{FLM2} uses some special features in the construction of the
moonshine module, in particular, ``triality'' (\cite{FLM1},
\cite{FLM2}). In fact, there is a conceptual way to define the vertex
operator map which is motivated by the $S_{3}$-symmetry of the Jacobi
identity and contragredient modules and which works also in much more
general cases (see \cite{FHL} and also \cite{DGM} in physicists'
language).  The hard part is to prove that the moonshine module
together with this abstractly defined vertex operator map is a vertex
operator algebra. This was first proved directly (i.e., without using
triality, as had been done in \cite{FLM2}) in \cite{DGM} using
techniques developed in string theory. Recently, this has also
been proved conceptually by the author \cite{H7} using the tensor
product theory for modules for a vertex operator algebra developed by
Lepowsky and the author \cite{HL1} \cite{HL4}--\cite{HL6} \cite{H6}
\cite{H7.5}
and some results of Dong-Mason-Zhu \cite{DMZ} and
Dong \cite{D} on modules for the vertex operator
algebra $V_{\Lambda}^{+}$. (Note that in more general cases in which
we can still define the vertex operator maps abstractly, it is not
always true that we will obtain a vertex operator algebra. In
\cite{DGM}, a vertex operator algebra is obtained in a family of cases
generalizing \cite{FLM2}.)

We now turn to the main subjects of this section.  At the end of this section,
we shall apply these results  to the problem above. The meterial below
in Section 1 is essentially taken from \cite{FHL}.

We first list some properties of the formal $\delta$-function
and
some easy consequences of the definition of vertex
operator algebra. First there is the fundamental property of the
$\delta$-function:
\begin{equation}\label{1-1}
f(x)\delta (x) = f(1)\delta (x) \;\;\mbox{for}\;\;f(x) \in
\Bbb{C}[x,x^{-1}].
\end{equation}
This property has many variants; in general, whenever an expression is
multiplied by the $\delta $-function, we may formally set the argument
appearing in the $\delta $-function equal to 1, provided the relevant algebraic
expressions make sense.
There are two basic identities for the $\delta$-function:
\begin{eqnarray}\label{1-2}
x^{-1}_1\delta \left( {x_2+x_0\over x_1}\right)  =
x^{-1}_2\delta \left( {x_1-x_0\over x_2}\right),
\end{eqnarray}
\begin{eqnarray}\label{1-3}
x^{-1}_0\delta \left( {x_1-x_2\over z_0}\right)  -
x^{-1}_0\delta \left( {x_2-x_1\over -x_0}\right)
=
x^{-1}_2\delta \left( {x_1-x_0\over x_2}\right).
\end{eqnarray}
Let $(V, Y, \bold{1}, \omega)$ be a vertex operator algebra.
We have the following immediate consequences of
the definition of vertex operator algebra:
\begin{eqnarray}
{[L(-1),Y(v,x)]} &=& Y(L(-1)v,x),\label{1-4}\\
{[L(0), Y(v, x)]} &=& Y(L(0)v,x) + xY(L(-1)v,x),\label{1-5}\\
{[L(1),Y(v,x)]} &=& Y(L(1)v,x) + 2xY(L(0)v,x) + x^2Y(L(-1)v,x)\label{1-6}
\end{eqnarray}
for any $v\in V$.
{}From the  $L(-1)$-derivative property and bracket
formulas (\ref{1-4}), we obtain
\begin{equation}\label{1-7}
e^{x_0L(-1)}Y(v,x)e^{-x_0L(-1)} = Y(e^{x_0L(-1)}v,x) = Y(v,x+x_0)
\end{equation}
Applying (\ref{1-7}) to $\bold{1}$ and then taking the constant term in
$x_{0}$,
we have
\begin{equation}\label{1-8}
Y(v,x)\bold{1} = e^{xL(-1)}v.
\end{equation}
Finally, one very important consequence is
 the skew-symmetry, that is, for any $u, v\in V$,
\begin{equation}\label{1-9}
Y(u, x)v=e^{xL(-1)}Y(v, -x)u.
\end{equation}
We derive (\ref{1-9}) as follows: We have
\begin{eqnarray}\label{1-10}
\lefteqn{x_0^{-1}\delta \left( {x_1-x_2\over x_0}\right)
 Y(u,x_1)Y(v,x_2)}\nno\\
 &&\quad-
x_0^{-1}\delta \left( {x_2-x_1\over -x_0}\right) Y(v,x_1)Y(u,x_2)\nno\\
&&=(-x_0)^{-1}\delta \left( {x_2-x_1\over -x_0}\right)
 Y(v,x_2)Y(u,x_1)\nno\\
&& \quad-
(-x_0)^{-1}\delta \left( {x_1-x_2\over -(-x_0)}\right) Y(u,x_1)Y(v,x_2).
\end{eqnarray}
By the Jacobi identity and (\ref{1-10}),
\begin{equation}
 x^{-1}_2\delta \left( {x_1-x_0\over x_2}\right) Y(Y(u,x_0)v,x_2)
=x^{-1}_1\delta \left( {x_2-(-x_0)\over x_1}\right) Y(Y(v,-x_0)u,x_1).
\end{equation}
Using the fundamental property of the $\delta$-function and the
identity (\ref{1-2}),
we obtain
\begin{equation}\label{1-13}
 x^{-1}_2\delta \left( {x_1-x_0\over x_2}\right) Y(Y(u,x_0)v,x_2)
= x^{-1}_2\delta \left( {x_1-x_0\over x_2}\right)  Y(Y(v,-x_0)u,x_2+x_{0}).
\end{equation}
In particular (taking the coefficient of $x_{1}^{-1}$ in (\ref{1-13})),
\begin{equation}\label{1-14}
Y(Y(u,x_0)v,x_2)=Y(Y(v,-x_0)u,x_2+x_{0}).
\end{equation}
But by the second equality in (\ref{1-7}),
\begin{equation}\label{1-15}
Y(Y(v,-x_0)u,x_2+x_{0})=Y(e^{x_{0}L(-1)}Y(v, -x_{0})u, x_{2}).
\end{equation}
By the creation property, (\ref{1-14}) and (\ref{1-15}),
\begin{eqnarray}
Y(u, x_{0})v&=&\lim_{x_{2}\to 0}Y(Y(u,x_0)v,x_2)\bold{1}\nno\\
&=&\lim_{x_{2}\to 0}Y(e^{x_{0}L(-1)}Y(v, -x_{0})u, x_{2})\bold{1}\nno\\
&=&e^{x_{0}L(-1)}Y(v, -x_{0})u.
\end{eqnarray}

Now we discuss the $S_{3}$-symmetry of the Jacobi
identity.  For the Jacobi identity for Lie algebras, if we call
\begin{equation}
[u, [v, w]]-[v, [u, w]]=[[u, v], w]
\end{equation}
``the Jacobi identity for the ordered triple $(u, v, w)$,''
then the Jacobi identity for $(u, v, w)$ implies the Jacobi identity for
any permutation of the ordered triple $(u, v, w)$. The $S_{3}$-symmetry for
the Jacobi identity for vertex operator algebra is an analogous statement.
(The analogy between Lie algebras and vertex operator algebras is the
reason why Frenkel, Lepowsky and Meurman called the main axiom for vertex
operator algebras the ``Jacobi identity.'' It would be more accurate
and less confusing to call this
identity the Frenkel-Lepowsky-Meurman identity or simply the FLM identity.)
Let us retain the axioms for a vertex operator algebra
except for the Jacobi identity, and let us call
\begin{eqnarray}\label{1-18}
&{\dps x_0^{-1}\delta
\left( \frac{x_1-x_2}{x_0}\right) Y(u,x_1)Y(v,x_2)w -
x_0^{-1}\delta \left( \frac{x_2-x_1}{-x_0}\right)
Y(v,x_2)Y(u,x_1)w}&\nno\\
&={\dps  x^{-1}_2\delta \left( \frac{x_1-x_0}{x_2}\right) Y(Y(u,x_0)v,x_2)w}&
\end{eqnarray}
``the Jacobi identity for the ordered triple $(u,v,w)."$ We also
assume that the consequences (\ref{1-7}) and (\ref{1-9}) hold.
By skew-symmetry (\ref{1-9})
for the pair $(u,v)$ and the second equality in (\ref{1-7}) for the vector
$Y(v,-x_0)u$ we have
\begin{equation}\label{1-19}
Y(Y(u,x_0)v, x_2) =
Y(e^{x_0L(-1)}Y(v,-x_0)u,x_2) = Y(Y(v,-x_0)u,x_2+x_0).
\end{equation}
Thus {}from (\ref{1-19}) and the identity  (\ref{1-2}),
the Jacobi identity (\ref{1-18}) for $(u,v,w)$ gives
\begin{eqnarray}
&{\dps (-x_0)^{-1}\delta
\left( \frac{x_2-x_1}{-x_0}\right) Y(v,x_2)Y(u,x_1)w -
(-x_0)^{-1}\delta \left( \frac{x_1-x_2}{-(-x_0)}\right)
Y(u,x_1)Y(v,x_2)}w& \nno\\
&= {\dps x^{-1}_1\delta \left(
\frac{x_2-(-x_0)}{x_1}\right) Y(Y(v,-x_0)u,x_1)w,}&
\end{eqnarray}
which is the Jacobi identity for  $(v,u,w)$  (with  $(x_1,x_2,x_0)$
replaced by  $(x_2,x_1, -x_0)$).

On the other hand, multiplying both sides of the Jacobi identity (\ref{1-18})
for  $(u,v,w)$
by  $e^{-x_2L(-1)}$ and using (\ref{1-9}) for the pairs  $(v,w)$,
$(v,Y(u,x_1)w)$
and  $(Y(u,x_0)v,w)$  and the outer equality in (\ref{1-7}) for
the vector  $u$,
we obtain
\begin{eqnarray}\label{1-21}
&{\dps x^{-1}_0\delta \left(
\frac{x_1-x_2}{x_0}\right)Y(u,x_1-x_2)Y(w,-x_2)v - x^{-1}_0\delta
\left(\frac{x_2-x_1}{-x_0}\right) Y(Y(u,x_1)w,-x_2)v}&\nno\\
&={\dps
x^{-1}_2\delta \left( \frac{x_1-x_0}{x_2}\right) Y(w,-x_2)Y(u,x_0)v.}&
\end{eqnarray}
Using the fundamental property of the $\delta$-function and (\ref{1-2}), we
can write  (\ref{1-21}) as
\begin{eqnarray}
&{\dps x^{-1}_1\delta \left( \frac{x_0+x_2}{x_1}\right)
Y(u,x_0)Y(w,-x_2)v + x^{-1}_2\delta \left( \frac{x_0-x_1}{-x_2}\right)
Y(Y(u,x_1)w,-x_2)v}& \nno\\
&={\dps x^{-1}_1\delta \left(
\frac{-x_2-x_0}{-x_1}\right) Y(w,-x_2)Y(u,x_0)v,}&
\end{eqnarray}
that is,
\begin{eqnarray}
\lefteqn{x^{-1}_1\delta \left( \frac{x_0-(-x_2)}{x_1}\right)
Y(u,x_0)Y(w,-x_2)v}
\nno \\
&& - x^{-1}_1\delta
\left( \frac{(-x_2)-x_0}{-x_1}\right) Y(w,-x_2)Y(u,x_0)v\nno \\
&&= (-x_2)^{-1}\delta \left( \frac{x_0-x_1}{-x_2}\right) Y(Y(u,x_1)w,-x_2)v,
\end{eqnarray}
the Jacobi identity for  $(u,w,v)$  (and  $(x_0,-x_2,x_1))$.  Since the two
permutation above of $(u, v, w)$ generate $S_{3}$, the permutation group
of $(u, v, w)$, we conclude:
\begin{propo}
 Under the assumptions indicated in the argument
above, the Jacobi identity for an ordered triple implies the
Jacobi identity for any
permutation of this triple.
\end{propo}

We turn next to the contragredient module for a module for a
vertex operator algebra.
Let  $(W,Y)$,  with
\begin{equation}
W = \coprod_{n\in \Bbb{C}}W_{(n)},
\end{equation}
be a module for a vertex operator algebra  $(V,Y, \bold{1},\omega),$
\begin{equation}
W' = \coprod_{n\in \Bbb{C}}W^*_{(n)}
\end{equation}
the graded dual space of  $W$ and $\langle \cdot, \cdot\rangle$ the pairing
between $W'$ and $W$.  We define the
{\it  contragredient  vertex operators}  $Y'(v,x)$ ($v \in  V$)
by means of the linear map
\begin{eqnarray}
V &\rightarrow & (\text{End}\ W')[[x,x^{-1}]]\nno\\
v &\mapsto & Y'(v,x) = \sum_{n\in
\Bbb{Z}}v'_nx^{-n-1} \;\;\;(\mbox{where}  \;\;v'_n \in
\mbox{\rm End}\ W'),
\end{eqnarray}
determined by the condition
\begin{equation}\label{1-27}
\langle Y'(v,x)w',w\rangle  =
\langle w',Y(e^{xL(1)}(-x^{-2})^{L(0)}v,x^{-1})w\rangle
\end{equation}
for  $v \in  V$,  $w' \in  W'$,  $w \in  W$.   The operator
$(-x^{-2})^{L(0)}$ has the obvious meaning; it acts on a vector of weight  $n
\in  \Bbb{Z}$  as multiplication by  $(-x^{-2})^n$.  Also note that
$e^{xL(1)}(-x^{-2})^{L(0)}v$  involves only finitely many (integral) powers of
$z$,  that the right-hand side of (\ref{1-27}) is a Laurent polynomial in  $x$,
 and
that the components  $v'_n$  of the formal Laurent series
$Y'(v,x)$  defined by (\ref{1-27}) indeed preserve  $W'.$

We give the space  $W'$ a  $\Bbb{C}$-grading by setting
\begin{equation}
W'_{(n)} = W^*_{(n)}\;\; \mbox{\rm for}\;\;  n \in  \Bbb{C}.
\end{equation}
The following proposition defines the  $V$-{\it module  contragredient to
$W$}:
\begin{theo}
The pair  $(W',Y')$  carries the structure of a
$V$-module.
\end{theo}
\begin{pf}
The axioms on the
grading
are clear.
For the Virasoro algebra properties, we note that
\begin{equation}
\langle Y'(\omega ,x)w',w\rangle  =
\langle w',Y(x^{-4}\omega ,x^{-1})w\rangle
\end{equation}
since
\begin{equation}
L(1)\omega  = L(-1)L(-2)\bold{1}=L(-2)L(-1)\bold{1}=0.
\end{equation}
Thus, defining component operators  $L'(n)$  by
\begin{equation}
Y'(\omega ,x) = \sum_{n\in
\Bbb{Z}}L'(n)x^{-n-2},
\end{equation}
we have
\begin{eqnarray}
\lefteqn{\langle \sum_{n\in
\Bbb{Z}}L'(n)x^{-n}w',w\rangle=}\nno \\
&& = \langle x^2Y'
(\omega ,x)w',w\rangle\nno \\
&&= \langle w',x^{-2}Y(\omega ,x^{-1})w\rangle\nno\\
&&= \langle w', \sum_{n\in \Bbb{Z}}L(-n)x^{-n}w\rangle,
\end{eqnarray}
and so
\begin{equation}
\langle L'(n)w',w\rangle  = \langle w',L(-n)w\rangle \;\;
\mbox{\rm for}\;\;
 n \in  \Bbb{Z}.
\end{equation}
This immediately gives us the Virasoro commutator relation for
$L'(n)$, $n\in \Bbb{Z}$.

We shall give  proofs of the
Jacobi identity and the  $L(-1)$-derivative property.
For these two axioms, we shall use some
commutator formulas motivated by the Lie group  SL$(2,\Bbb{C})$,  but
formulated and proved in terms of formal series. We shall omit the
 proofs of these formulas;
they can be found in \cite{FHL} and are all direct
calculations.

\begin{lemma}  Let
\begin{equation}
f(x) \in  x\Bbb{C}[[x]].
\end{equation}
We have the following identities, valid on any module for the Lie algebra
$\frak{s}\frak{l}(2)$  spanned by  $L(-1)$, $L(0)$, $L(1):$
\begin{equation}
L(-1)e^{f(x)L(0)} = e^{f(x)L(0)}L(-1)e^{-f(x)},
\end{equation}
\begin{equation}
L(1)e^{f(x)L(0)} = e^{f(x)L(0)}L(1)e^{f(x)},
\end{equation}
\begin{eqnarray} \label{1-37}
\lefteqn{L(-1)e^{f(x)L(1)}=}\nno\\
&&= e^{f(x)L(1)}L(-1) - 2f(x)L(0)e^{f(x)L(1)} -
f(x)^2L(1)e^{f(x)L(1)}\nno\\
&&= e^{f(x)L(1)}L(-1) - 2f(x)e^{f(x)L(1)}L(0) + f(x)^2e^{f(x)L(1)}L(1).
\end{eqnarray}
These identities also hold for more general  $f$  for which the series are well
defined, such as
\begin{equation}
f(x,x_0) \in  x\Bbb{C}[[x,x_0]].
\end{equation}
\end{lemma}

Now we establish the  $L(-1)$-derivative property.  For convenience,
we assume that  $v \in  V$  is homogeneous of weight  $n \in  \Bbb{Z}:$
$L(0)v = nv.$
Using the definition  $Y'(\cdot ,x)$  and the chain rule we get
\begin{eqnarray} \label{1-39}
\lefteqn{\langle {d\over dx}Y'(v,x)w',w\rangle=}\nno\\
&&={d\over dx}\langle w',Y(e^{xL(1)}(-x^{-2})^{L(0)}v,x^{-1})w\rangle \nno\\
&&= \langle w',{d\over dx}Y(e^{xL(1)}(-x^{-2})^{L(0)}v,x^{-1})w\rangle \nno\\
&&= \langle w',Y({d\over dx}(e^{xL(1)}(-x^{-2})^{L(0)})v,x^{-1})w\rangle\nno \\
&&\;\;\;\;+ \langle w',{d\over dx}Y(v_1,x^{-1})\left| _{v_1 =
e^{xL(1)}(-x^{-2})^{L(0)}v}w\rangle,\right.
\end{eqnarray}
where  $w'$  and  $w$  are arbitrary elements of  $W'$ and  $W$,
respectively.  We perform the indicated calculations:
\begin{eqnarray} \label{1-40}
\lefteqn{{d\over dx}(e^{xL(1)}(-x^{-2})^{L(0)})=}\nno\\
&&= L(1)e^{xL(1)}(-x^{-2})^{L(0)} - 2x^{-1}e^{xL(1)}L(0)(-x^{-2})^{L(0)},
\end{eqnarray}
\begin{eqnarray} \label{1-41}
\lefteqn{{d\over dx}Y(v_1,x^{-1})\left| _{v_1 = e^{xL(1)}(-x^{-2})^{L(0)}v}
\right.=}\nno \\
&&= -x^{-2}{d\over dx^{-1}}Y(v_1,x^{-1})\left| _{v_1 =
e^{xL(1)}(-x^{-2})^{L(0)}v}\right.\nno \\
&&= -x^{-2}Y(L(-1)v_1,x^{-1})\left| _{v_1 = e^{xL(1)}(-x^{-2})^{L(0)}v}\right.
\nno\\
&&= -x^{-2}Y(L(-1)e^{xL(1)}(-x^{-2})^{L(0)}v,x^{-1})\nno\\
&&= -x^{-2}Y((e^{xL(1)}L(-1) - 2xe^{xL(1)}L(0)\nno\\
&&\;\;\;\;+ x^2L(1)e^{xL(1)})(-x^{-2})^nv,x^{-1})\nno\\
&&= Y(e^{xL(1)}(-x^{-2})^{n+1}L(-1)v,x^{-1})\nno\\
&&\;\;\;\;+ Y(2x^{-1}e^{xL(1)}L(0)(-x^{-2})^nv,x^{-1})\nno\\
&&\;\;\;\;- Y(L(1)e^{xL(1)}(-x^{-2})^nv,x^{-1})\nno\\
&&= Y(e^{xL(1)}(-x^{-2})^{L(0)}L(-1)v,x^{-1})\nno\\
&&\;\;\;\;+ Y(2x^{-1}e^{xL(1)}L(0)(-x^{-2})^{L(0)}v,x^{-1})\nno\\
&&\;\;\;\;- Y(L(1)e^{xL(1)}(-x^{-2})^{L(0)}v,x^{-1}).
\end{eqnarray}
Here we have used the outer equality in (\ref{1-37}) and
the fact that
\begin{equation}
L(0)L(-1)v = L(-1)(L(0) + 1)v = (n+1)L(-1)v.
\end{equation}
Substituting (\ref{1-40}) and (\ref{1-41}) into (\ref{1-39}) we get
\begin{eqnarray}
\lefteqn{\langle {d\over dx}Y'(v,x)w',w\rangle =}\nno\\
&&= \langle w',Y(L(1)e^{xL(1)}(-x^{-2})^{L(0)}v\nno\\
&&\;\;\;\;\;\;\;\;- 2x^{-1}e^{xL(1)}L(0)(-x^{-2})^{L(0)}v,x^{-1})w\rangle
\nno\\
&&\;\;\;\;+ \langle w',Y(e^{xL(1)}(-x^{-2})^{L(0)}L(-1)v,x^{-1})w\rangle \nno\\
&&\;\;\;\;+ \langle w',Y(2x^{-1}e^{xL(1)}L(0)(-x^{-2})^{L(0)}v,x^{-1})w\rangle
\nno\\
&&\;\;\;\;- \langle w',Y(L(1)e^{xL(1)}(-x^{-2})^{L(0)}v,x^{-1})w\rangle \nno\\
&&= \langle w',Y(e^{xL(1)}(-x^{-2})^{L(0)}L(-1)v,x^{-1})w\rangle \nno\\
&&= \langle Y'(L(-1)v,x)w',w\rangle,
\end{eqnarray}
proving the $L(-1)$-derivative property.

Finally, we shall prove the Jacobi identity.  Let  $v_1,v_2 \in  V$,
$w \in  W$  and  $w' \in  W'$.  What we want to prove can be written as
follows:
\begin{eqnarray} \label{1-44}
\lefteqn{\langle x^{-1}_0\delta \left( {x_1-x_2\over x_0}\right)
Y'(v_1,x_1)Y'(v_%
2,x_2)w',w\rangle }\nno\\
&&\;\;\;\;-
\langle x^{-1}_0\delta \left( {x_2-x_1\over -x_0}\right)
Y'(v_2,x_2)Y'(v_%
1,x_1)w',w\rangle\nno\\
&&=
\langle x^{-1}_2\delta \left( {x_1-x_0\over x_2}\right)
Y'(Y(v_1,x_0)v_2,x_2%
)w',w\rangle.
\end{eqnarray}
But by the definition (\ref{1-27}) of contragredient vertex operator, we have
\begin{eqnarray}
\lefteqn{\langle Y'(v_1,x_1)Y'(v_2,x_2)w',w\rangle }\nno \\
&&= \langle w',Y(e^{x_2L(1)}(-x^{-2}_2)^{L(0)}v_2,x^{-1}_2)
Y(e^{x_1L(1)}(-x^{-2}_1)
^{L(0)}v_1,x^{-1}_1)w\rangle
\end{eqnarray}
\begin{eqnarray}
\lefteqn{\langle Y'(v_2,x_2)Y'(v_1,x_1)w',w\rangle }\nno \\
&&=
\langle w',Y(e^{x_1L(1)}(-x^{-2}_1)^{L(0)}v_1,x^{-1}_1)Y(e^{x_2L(1)}(-x^{-2}_2)
^{L(0)}v_2,x^{-1}_2)w\rangle
\end{eqnarray}
\begin{eqnarray}
\lefteqn{\langle Y'(Y(v_1,x_0)v_2,x_2)w',w\rangle }\nno \\
&&= \langle w',Y(e^{x_2L(1)}(-x^{-2}_2)^{L(0)}Y(v_1,x_0)v_2,x^{-1}_2)w\rangle,
\end{eqnarray}
and {}from the Jacobi identity for  $W$  we have
\begin{eqnarray}
\lefteqn{\langle w',\left( {-x_0\over x_1x_2}\right) ^{-1}\delta
 \left( {x^{-1}_1-x^{-1}_2\over -x_0/x_1x_2} \right)
Y(e^{x_1L(1)}(-x^{-2}_1)^{L(0)}v_1,x^{-1}_1)\cdot}\nno \\
&&\;\;\;\;\;\;\;\;\;\;\cdot
Y(e^{x_2L(1)}(-x^{-2}_2)^{L(0)}v_2,x^{-1}_2)w\rangle\nno \\
&&\;\;\;\;- \langle w',\left( {-x_0\over x_1x_2}\right) ^{-1}\delta
\left( {x^{-1}_2-x^{-1}_1\over x_0/x_1x_2}\right)
Y(e^{x_2L(1)}(-x^{-2}_2)^{L(0)}v_2,x^{-1}_2)\cdot \nno\\
&&\;\;\;\;\;\;\;\;\;\;\cdot
Y(e^{x_1L(1)}(-x^{-2}_1)^{L(0)}v_1,x^{-1}_1)w\rangle \nno\\
&&= \langle w',
(x^{-1}_2)^{-1}\delta \left( {x^{-1}_1+x_0/x_1x_2\over x^{-1}_2}%
\right) \cdot\nno \\
&&\;\;\;\;\;\;\;\;\;\;\cdot
Y(Y(e^{x_1L(1)}(-x^{-2}_1)^{L(0)}v_1,-x_0/x_1x_2)e^{x_2L(1)}(-x^{%
-2}_2)^{L(0)}v_2,x^{-1}_2)w\rangle,\nno\\
&&
\end{eqnarray}
or equivalently,
\begin{eqnarray}
\lefteqn{-
\langle w',x^{-1}_0\delta \left( {x_2-x_1\over -x_0}\right)
Y(e^{x_1L(1)}(-x^{-%
2}_1)^{L(0)}v_1,x^{-1}_1)\cdot}\nno \\
&&\;\;\;\;\;\;\;\;\;\;\cdot
Y(e^{x_2L(1)}(-x^{-2}_2)^{L(0)}v_2,x^{-1}_2)w\rangle\nno \\
&&\;\;\;\;+
\langle w',x^{-1}_0\delta
\left( {x_1-x_2\over x_0}\right)
Y(e^{x_2L(1)}(-x^{-2}_2)^{L(0)}v_2,x^{-1}_2)\cdot\nno \\
&&\;\;\;\;\;\;\;\;\;\;\cdot
Y(e^{x_1L(1)}(-x^{-2}_1)^{L(0)}v_1,x^{-1}_1)w\rangle\nno \\
&&= \langle w',x^{-1}_1\delta \left( {x_2+x_0\over x_1}\right) \cdot\nno\\
&&\;\;\;\;\;\;\;\;\;\;\cdot
Y(Y(e^{x_1L(1)}(-x^{-2}_1)^{L(0)}v_1,-x_0/x_1x_2)e^{x_2L(1)}(-x^{%
-2}_2)^{L(0)}v_2,x^{-1}_2)w\rangle.\nno\\
&&
\end{eqnarray}
(As usual, the reader should be observing that the formal Laurent series which
arise are well defined.)  Thus (by (\ref{1-2})) the desired result
(\ref{1-44}) is
equivalent to
\begin{eqnarray}
\lefteqn{x^{-1}_1\delta \left( {x_2+x_0\over x_1}\right)
Y(e^{x_2L(1)}(-x^{-2}_2)^{L(0)
}Y(v_1,x_0)v_2,x^{-1}_2) }\nno\\
&&=
x^{-1}_1\delta \left( {x_2+x_0\over x_1}\right)
Y(Y(e^{x_1L(1)}(-x^{-2}_1)^{L(0%
)}v_1,-x_0/x_1x_2)\cdot\nno \\
&&\;\;\;\;\;\;\;\;\;\;\cdot e^{x_2L(1)}(-x^{-2}_2)^{L(0)}v_2,x^{-1}_2),
\end{eqnarray}
or to
\begin{eqnarray}
\lefteqn{Y(e^{x_2L(1)}(-x^{-2}_2)^{L(0)}Y(v_1,x_0)v_2,x^{-1}_2)}\nno \\
&&=
Y(Y(e^{(x_2+x_0)L(1)}(-(x_2+x_0)^{-2})^{L(0)}v_1,-x_0/(x_2+x_0)x_2)%
\cdot\nno \\
&&\;\;\;\;\;\;\;\;\;\;\cdot e^{x_2L(1)}(-x^{-2}_2)^{L(0)}v_2,x^{-1}_2).
\end{eqnarray}
If we can prove
\begin{eqnarray}
\lefteqn{e^{x_2L(1)}(-x^{-2}_2)^{L(0)}Y(v_1,x_0)}\nno\\
&&=
Y(e^{(x_2+x_0)L(1)}(-(x_2+x_0)^{-2})^{L(0)}v_1,-x_0/(x_2+x_0)x_2)\cdot\nno\\
&&\;\;\;\;\;\;\;\;\;\;\cdot e^{x_2L(1)}(-x^{-2}_2)^{L(0)}
\end{eqnarray}
or equivalently, the conjugation formula
\begin{eqnarray}
\lefteqn{e^{xL(1)}(-x^{-2})^{L(0)}Y(v,x_0)(-x^{-2})^{-L(0)}e^{-xL(1)}}\nno \\
&&= Y(e^{(x+x_0)L(1)}(-(x+x_0)^{-2})^{L(0)}v,-x_0/(x+x_0)x)
\end{eqnarray}
for any element  $v$  of a vertex operator algebra, where the operators act on
the algebra itself, then we will be done.  But for this, it is sufficient to
prove the following lemma:
\begin{lemma} Let  $V$  be a vertex operator algebra.
 The following conjugation formulas hold on  $V:$
\begin{eqnarray}
x^{L(0)}Y(v,x_0)x^{-L(0)} &=& Y(x^{L(0)}v,xx_0)\\
e^{xL(1)}Y(v,x_0)e^{-xL(1)} &= &
Y(e^{x(1-xx_0)L(1)}(1-xx_0)^{-2L(0)}v,x_0/(1-xx_0)).
\end{eqnarray}
\end{lemma}
The proof of this lemma, which we omit here, can be found in \cite{FHL}.
This finishes the proof of the theorem.
\end{pf}

The functor taking a $V$-module to its contragredient module has some
important properties which we state without proof (see \cite{FHL}):

\begin{propo}
There is a natural isomorphism between the double contragredient module
$(W'', Y'')$ and $(W, Y)$.
\end{propo}

\begin{propo}
The module $(W, Y)$ is irreducible if and only if $(W', Y')$ is irreducible.
\end{propo}

\begin{propo}\label{biform}
The module $(W, Y)$ is isomorphic to its contragredient module $(W', Y')$ if
and only if there exists a nondegenerate bilinear form $(\cdot, \cdot)_{W}$
on $W$ such that
\begin{equation}\label{1-55}
(W_{(m)}, W_{(n)})_{W}=0,\;\;\;m\ne n
\end{equation}
and
\begin{equation}\label{1-56}
(Y'(v,x)w_{1}, w_{2})_{W}  =
 (w_{1}, Y(e^{xL(1)}(-x^{-2})^{L(0)}v,x^{-1})w_{2})_{W}.
\end{equation}
 If $V$
as a $V$-module is isomorphic to $V'$, the bilinear form $(\cdot, \cdot)_{V}$
is symmetric.
\end{propo}

We return to our problem of defining the vertex operator map for
$V^{\natural}$. Let $V$ be a vertex operator algebra and $W$ a $V$-module.
Assume that both $V$ and $W$ as $V$-modules are isomorphic to themselves.
By Proposition \ref{biform}, there are nondegenerate bilinear forms
$(\cdot, \cdot)_{V}$ and $(\cdot, \cdot)_{W}$ satisfying the two conditions
in Proposition \ref{biform}. In addition, $(\cdot, \cdot)_{V}$ is
symmetric. Assume that there is a vertex operator map
$$Y_{V\oplus W}: (V\oplus W)\otimes (V\oplus W)\to
(\text{End}\; (V\oplus W))[[x, x^{-1}]]$$
such that $(V\oplus W, Y_{V\oplus W}, \bold{1}, \omega)$ ($\bold{1}$ and
$\omega$ are the vacuum and the Virasoro element of $V$, respectively)
is a vertex operator algebra satisfying the following:

\begin{enumerate}

\item The vertex
operator algebra structure on $V$ and the module structure on $W$ are
substructures of it.

\item As a module for itself, it is isomorphic its contragredient
module and the corresponding symmetric nondegenerate bilinear form
$(\cdot, \cdot)_{V\oplus W}$ is
defined by
\begin{equation}
((v_{1}, w_{2}), (v_{2}, w_{2}))_{V\oplus W}=(v_{1}, v_{2})_{V}
+(w_{1}, w_{2})_{W}
\end{equation}
for all $v_{1}, v_{2}\in V$ and $w_{1}, w_{2}\in W$.

\item  The involution which is the identity on $V$ and is $-1$ on $W$ is an
automorphism of $(V\oplus W, Y_{V\oplus W}, \bold{1}, \omega)$.

\end{enumerate}

Then we must have the
following:

\begin{enumerate}

\item The module $W$ is $\Bbb{Z}$-graded.

\item The bilinear form $(\cdot, \cdot)_{W}$ is symmetric.

\item We have the following formulas: For any $v\in V$ and $w\in W$,
\begin{equation}\label{1-58}
Y_{V\oplus W}(w, x)v=e^{xL(-1)}Y_{W}(v, -x)w
\end{equation}
and
for any $v\in V$, $w_{1}, w_{2}, w_{3}\in W$,
\begin{eqnarray}
(w_{3}, Y_{V\oplus W}(w_{1}, x)w_{2})_{W}&=&0,\label{1-59}\\
(v, Y_{V\oplus W}(w_{1}, x)w_{2})_{V}&=&(Y_{W}(v, -x^{-1})
e^{xL(1)}(-x^{2})^{-L(0)}w_{1},
e^{x^{-1}L(1)}w_{2})_{W},\nno\\
&&\label{1-60}
\end{eqnarray}
where $Y_{V}$ and $Y_{W}$ are the vertex operator maps for $V$ and $W$,
respectively.
\end{enumerate}

We see that the vertex operator map $Y_{V\oplus W}$ is determined completely
by the vertex operator maps $Y_{V}$, $Y_{W}$, the bilinear forms
$(\cdot, \cdot)_{V}$, $(\cdot, \cdot)_{W}$  and (\ref{1-58})--(\ref{1-60}).
Thus even if we do not know whether $V\oplus W$ is such a vertex operator
algebra, we can still define a vertex operator map $Y_{V\oplus W}$ using
$Y_{V}$, $Y_{W}$, the bilinear forms
$(\cdot, \cdot)_{V}$, $(\cdot, \cdot)_{W}$
and (\ref{1-58})--(\ref{1-60}). In particular, since $V^{+}_{\Lambda}$ and
$(V_{\Lambda}^{T})^{+}$ as $V^{+}_{\Lambda}$-modules are both isomorphic
to their contragredient modules, we can define a vertex operator map
$Y_{V^{\natural}}$ for
$V^{\natural}=V^{+}_{\Lambda}\oplus (V_{\Lambda}^{T})^{+}$.

\section{Intertwining operators, fusion rules and Verlinde algebras}

We first define intertwining operators and fusion rules for a
vertex operator following \cite{FHL}.

Let
\begin{equation}
V\{x\} = \left\lbrace \sum_{n
\in
\Bbb{C}}v_nx^n|v_n \in  V\right\rbrace
\end{equation}
be the vector space of  $V$-valued  formal series involving
the complex powers  of  $x$ with coefficients in  a
vector space  $V.$
\begin{defi}
Let  $V$  be a vertex operator algebra and let
$(W_1,Y_1)$,  $(W_2,Y_2)$  and  $(W_3,Y_3)$  be three  $V$-modules  (not
necessarily distinct, and possibly equal to  $V)$.  An {\it intertwining
operator of type  ${3}\choose{12}$}  (or {\it of type
${W_3}\choose{W_1\ W_2} $})  is a linear map  $W_1\otimes W_2
\rightarrow  W_3\{x\}$,  or equivalently,
\begin{eqnarray}
W_1 &\rightarrow & (\mbox{\rm Hom}(W_2,W_3))\{x\}\nno\\
w & \mapsto & \cal{Y}(w,x) =\sum_{n\in
\Bbb{Q}}w_nx^{-n-1}\;\;\; (\mbox{\rm where}\;\;
w_n \in  \mbox{\rm Hom}(W_2,W_3))
\end{eqnarray}
such that ``all the defining properties of a module action that make
sense hold"
(cf. the definition of $V$-module).
That is, for  $v \in  V$,  $w_{(1)} \in  W_1$ and
$w_{(2)} \in  W_2,$
\begin{equation}\label{2-3}
w_{(1)n}w_{(2)} = 0\;\;  \mbox{\rm for}\;\; n \;\;
\mbox{\rm whose real part is sufficiently large;}
\end{equation}
the following Jacobi identity holds for the operators  $Y_{i}(v,\cdot )$,
$i=1, 2, 3$,
$\cal{Y}(w_{(1)},\cdot )$  acting on the element  $w_{(2)}:$
\begin{eqnarray}
\lefteqn{x^{-1}_0\delta \left( {x_1-x_2\over x_0}\right)
Y_3(v,x_1)\cal{Y}(w_{(1)},x_2%
)w_{(2)}}\nno\\
&&\;\;\;\;-
x^{-1}_0\delta \left( {x_2-x_1\over -x_0}\right)
\cal{Y}(w_{(1)},x_2)Y_2(v,x_1%
)w_{(2)} \nno\\
&&=
x^{-1}_2\delta \left( {x_1-x_0\over x_2}\right) \cal{Y}(Y_1(v,x_0)w_{(1)},x_2)%
w_{(2)}
\end{eqnarray}
(note that the first term on the left-hand side is algebraically meaningful
because of condition (\ref{2-3}), and the other terms are meaningful by the
usual
properties of modules; also note that this Jacobi identity involves integral
powers of  $x_0$ and  $x_1$ and complex powers of  $x_2);$
\begin{equation}
{d\over dx}\cal{Y}(w_{(1)},x) = \cal{Y}(L(-1)w_{(1)},x),
\end{equation}
where  $L(-1)$  is the operator acting on  $W_{1}.$
\end{defi}

We may denote the intertwining operator just defined by
\begin{equation}
\cal{Y}^3_{12}\;\;\; \mbox{\rm or}\;\;\cal{Y}^{W_3}_{W_1W_2},
\end{equation}
if necessary, to indicate its type.

Note that  $Y(\cdot ,x)$  acting on  $V$
is an example of
an intertwining operator of type  $ V\choose {V\ V} $,  and
$Y(\cdot ,x)$  acting on a  $V$-module  $W$  is an example of an intertwining
operator of type  $ W\choose {V\ W} $.  These intertwining operators
satisfy the normalization condition  $Y(\bold{1},x) = 1$.

The intertwining operators of type  $ 3\choose {1\ 2}$  clearly form a
vector space, which we denote by  $\cal{V}^3_{12}$ or
$\cal{V}^{W_3}_{W_1W_2}$.  We set
\begin{equation}
N^3_{12} = N^{W_3}_{W_1W_2} = \dim \ \cal{V}^3_{12}
\;\;(\le \infty ).
\end{equation}
These numbers are called the
{\it fusion rules} associated with the algebra and
modules.  Then for example, assuming that  $V$  and the  $V$-module  $W$  are
nonzero, the corresponding fusion rules are positive:
\begin{eqnarray}
N^V_{VV} &\ge& 1,\\
N^W_{VW} &\ge&  1,\\
N^W_{WV} &\ge&  1.
\end{eqnarray}
In \cite{FHL} and \cite{HL5}, it is shown that the fusion rules have
the following symmetry property: Define
\begin{equation}
N_{ijk}=N_{W_{i}W_{j}W_{k}}=N^{W'_{k}}_{W_{i}W_{j}}
\end{equation}
for $i, j, k=1, 2, 3$.
Then
for any element $\sigma\in S_{3}$,
we have
\begin{equation}\label{2-10}
N_{\sigma(1)\sigma(2)\sigma(3)}=N_{123}.
\end{equation}

If the vertex operator algebra $V$ is rational, that is, V satisfies the
conditions:
(i) there are only finitely many irreducible $V$-modules
(up to equivalence), (ii) every $V$-module is completely
reducible, (iii) all the fusion rules are finite, then we can
define an algebra called the {\it fusion algebra} or the {\it Verlinde
algebra} using fusion rules for the irreducible modules as follows:
Assume that there are $m$ inequivalent irreducible $V$-modules. Let
$A$ be the abelian group tensor product of
the $K$-group  of the $V$-modules with
$\Bbb{C}$. Then $A$ has a natural structure of a vector space. Since $V$ is
rational, we have
\begin{equation}
A=\sum_{i=1}^{m}\Bbb{C}\phi_{i}
\end{equation}
where $\phi_{i}$, $i=1, \dots, m$, are all the equivalence classes containing
irreducible modules. We define a product on $A$ by
\begin{equation}
\phi_{i}\cdot \phi_{j}=\sum_{k=1}^{m}N_{ij}^{k}\phi_{k}
\end{equation}
for all $i, j=1, \dots, m$, where $N^{k}_{ij}$, $1\le i, j,k \le m$, are
the fusion rules $N^{W_{k}}_{W_{i}W_{j}}$ for any $W_{i}\in \phi_{i}$,
$W_{j}\in \phi_{j}$ and $W_{k}\in \phi_{k}$.
By the symmetry (\ref{2-10}), it is clear that this product is commutative.
When the intertwining operators for
the  vertex operator algebra satisfy certain additional
conditions, it can be proved
that this product is also associative. One condition that
we need is that all irreducible $V$-modules are $\Bbb{R}$-graded.
If $V$ is rational, then this condition implies that every $V$-module is
$\Bbb{R}$-graded, that is, the weight of an element of a $V$-module is
always a real number.
We also need
an additional condition.
Given any $V$-modules $W_{1}$, $W_{2}$, $W_{3}$, $W_{4}$ and $W_{5}$,
let $\cal{Y}_{1}$, $\cal{Y}_{2}$, $\cal{Y}_{3}$ and $\cal{Y}_{4}$
be intertwining operators of type ${W_{4}}\choose {W_{1}W_{5}}$,
${W_{5}}\choose {W_{2}W_{3}}$, ${W_{5}}\choose {W_{1}W_{2}}$ and
${W_{4}}\choose {W_{5}W_{3}}$, respectively. Consider the following
conditions for the product of $\cal{Y}_{1}$ and $\cal{Y}_{2}$ and
for the iterate of $\cal{Y}_{3}$ and $\cal{Y}_{4}$, respectively:

\begin{description}

\item[Convergence and extension property for products]
There exists
an integer $N$
(depending only on $\cal{Y}_{1}$ and $\cal{Y}_{2}$), and
for any $w_{(1)}\in W_{1}$,
$w_{(2)}\in W_{2}$, $w_{(3)}\in W_{3}$, $w'_{(4)}\in W'_{4}$, there exist
$j\in \Bbb{N}$, $r_{i}, s_{i}\in \Bbb{R}$, $i=1, \dots, j$, and analytic
functions $f_{i}(z)$ on $|z|<1$, $i=1, \dots, j$,
satisfying
\begin{equation}
\wt w_{(1)}+\wt w_{(2)}+s_{i}>N,\;\;\;i=1, \dots, j,
\end{equation}
such that
\begin{equation}
\langle w'_{(4)}, \cal{Y}_{1}(w_{(1)}, x_{2})
\cal{Y}_{2}(w_{(2)}, x_{2})w_{(3)}\rangle_{W_{4}}
\lbar_{x_{1}^{n}= e^{n\log z_{1}}, \;x_{2}^{n}=e^{n\log z_{2}},\;
n\in \Bbb{C}}
\end{equation}
is convergent when $|z_{1}|>|z_{2}|>0$ and can be analytically extended to
the multi-valued analytic function
\begin{equation}\label{phyper}
\sum_{i=1}^{j}z_{2}^{r_{i}}(z_{1}-z_{2})^{s_{i}}
f_{i}\left(\frac{z_{1}-z_{2}}{z_{2}}\right)
\end{equation}
when $|z_{2}|>|z_{1}-z_{2}|>0$.

\item[Convergence and extension property for iterates]
\hspace{.5em}There exists an integer
$\tilde{N}$
(depending only on $\cal{Y}_{3}$ and $\cal{Y}_{4}$), and
for any $w_{(1)}\in W_{1}$,
$w_{(2)}\in W_{2}$, $w_{(3)}\in W_{3}$, $w'_{(4)}\in W'_{4}$, there exist
$k\in \Bbb{N}$, $\tilde{r}_{i}, \tilde{s}_{i}\in \Bbb{R}$, $i=1, \dots, k$,
and analytic
functions $\tilde{f}_{i}(z)$ on $|z|<1$, $i=1, \dots, k$,
satisfying
\begin{equation}
\wt w_{(2)}+\wt w_{(3)}+\tilde{s}_{i}>\tilde{N},\;\;\;i=1, \dots, k,
\end{equation}
 such that
\begin{equation}
\langle w'_{(4)},
\cal{Y}_{4}(\cal{Y}_{3}(w_{(1)}, x_{0})w_{(2)}, x_{2})w_{(3)}\rangle_{W_{4}}
\lbar_{x^{n}_{0}=e^{n\log (z_{1}-z_{2})},\;x^{n}_{2}=e^{n\log z_{2}},\;
n\in \Bbb{C}}
\end{equation}
is convergent when $|z_{2}|>|z_{1}-z_{2}|>0$ and can be analytically extended
to the multi-valued analytic function
\begin{equation}
\sum_{i=1}^{k}z_{1}^{\tilde{r}_{i}}z_{2}^{\tilde{s}_{i}}
\tilde{f}_{i}\left(\frac{z_{2}}{z_{1}}\right)
\end{equation}
when $|z_{1}|>|z_{2}|>0$.

\end{description}

If for any $V$-modules $W_{1}$, $W_{2}$, $W_{3}$, $W_{4}$ and $W_{5}$ and
 any intertwining operators $\cal{Y}_{1}$ and $\cal{Y}_{2}$
of the types as above, the convergence and extension property of products
holds,
we say that
{\it  the products of the
 intertwining operators for $V$ have the convergence and extension property}.
Similarly we can define the meaning of the phrase
 {\it the iterates of the intertwining
operators for $V$ have the convergence and extension property}.

We also need the notion of generalized module: A {\it generalized $V$ module}
is
a pair $(W, Y)$ satisfying all the axioms for a $V$-module except
the two grading axioms: $\dim W_{(n)}< \infty$ for all $n\in \Bbb{C}$ and
$W_{(n)}=0$ for $n\in \Bbb{C}$ whose real part is sufficiently small.
If a generalized $V$-module $W=\coprod_{n\in \Bbb{C}}W_{(n)}$
satisfies the second grading axiom above,
we say that $W$ is {\it lower-truncated}.
We have the following
result:
\begin{theo}
Let $V$ be a rational vertex operator algebra for which all irreducible modules
are $\Bbb{R}$-graded. Assume
that $V$ satisfies the following conditions:
\begin{enumerate}
\item Every finitely-generated lower-truncated generalized $V$-module
is a $V$-module.

\item The products or the iterates of
the intertwining operators for $V$ have the convergence and extension property.
\end{enumerate}
Then the Verlinde algebra for $V$ is a commutative associative algebra
with unit.
\end{theo}

This theorem is an easy consequence of the associativity of the tensor product
theory for modules for a vertex operator algebra developed by Lepowsky and
the author \cite{HL1} \cite{HL4}--\cite{HL6} \cite{H6}.

Fusion rules and Verlinde algebras are very important concepts and
tools in the study of conformal field theory. One of the most
interesting results
in the mathematical study of conformal field theory is that
the fusion rules and their higher-genus generalizations
 for the WZNW conformal field theory can be expressed
 in terms of elementary functions (actually, the sine functions)
\cite{V}.
On the other hand, these fusion rules and generalizations can also be
shown to be equal to the dimensions of the space of ``generalized
theta functions'' on the moduli spaces of semistable principal bundles
on smooth projective irreducible algebraic curves \cite{KNR}.  Thus
one obtains a simple and beautiful formula for these dimensions. These
are the so called Verlinde formulas.  Mathematical proofs of these
formulas have been obtained in \cite{TUY} and
\cite{F}.

\section{Geometric interpretation of vertex operator algebras}

We give a brief description of the geometric interpretation of vertex operator
algebras in this section. The geometric interpretations of vertex operators,
their duality properties and their transformation properties under the
 projective transformations were first given by Frenkel \cite{Fr} using
the geometry of $\Bbb{C}\cup \{\infty\}$ with some discs deleted. The
complete geometric interpretation is obtained in \cite{H1} and \cite{H8}.
The formulation using operads is given in \cite{HL2} and \cite{HL3}.
See \cite{H1}--\cite{H4}, \cite{H8}, \cite{HL2} and \cite{HL3} for
details and other expositions.

In classical algebraic theories we study mostly algebraic structures defined by
binary operations. These binary operations can always be described by
one-dimensional geometric objects. For example,
Lie algebras can be described by
binary trees. A Lie algebra can be defined to be
a ``linear representation" of the moduli space of binary trees with a ``welding
operation," satisfying certain ``conservation" and ``orientation" properties
\cite{H1} \cite{H5}.
Any
associative binary operation, for example, the multiplication for a group or
an algebra, can be described using
the moduli space of circles with punctures and
local coordinates \cite{HL2} \cite{HL3}.
The general philosophy behind the geometric interpretation of vertex
operator algebras is to study certain two-dimensional analogues of the
classical binary operations,
that is, to study operations described by two-dimensional analogues of binary
trees or circles with punctures and local coordinates.

The two-dimensional analogues, used to describe vertex operator
algebras, of both binary trees
and circles with punctures and local coordinates are spheres with
analytically parametrized boundaries,
where by spheres we mean one-dimensional compact connected genus-zero complex
manifolds. These spheres with boundaries are in some sense
equivalent to  spheres with ordered points (which are called punctures),
one negatively oriented and others positively oriented,
and local coordinates vanishing at these
points, as is explained in \cite{H1} and \cite{H8}.
We will use the the index $0$
to denote the negatively oriented puncture on such a sphere
with punctures and local coordinates. Let $S_{1}$ and $S_{2}$ be two
such spheres
with punctures and local coordinates,
$p_{j}$, $j=0, \dots, m$,  the punctures of $S_{1}$,
$q_{k}$, $k=0, \dots, n$, the
punctures of $S_{2}$, $(U_{j}, \varphi_{j})$, $j=0, \dots, m$,  the
local coordinates vanishing at $p_{j}$ and
$(V_{k}, \psi_{k})$, $k=0, \dots, n$,  the
local coordinates vanishing at $q_{k}$.
For any integer $i$ satisfying $0<i\le n$,
we would like to sew $S_{1}$ and $S_{2}$ through the $i$-th puncture of $S_{1}$
and the $0$-th puncture of $S_{2}$ to obtain
a new spheres with punctures and local coordinates.
Assume that there exists a
positive number $r$ such that $\varphi _{i}(U_{i})$ contains the
closed disc $\bar{B}_{0}^{r}$ centered at $0$ with radius $r$ and
$\psi_{0}(V_{0})$
contains the closed disc $\bar{B}_{0}^{1/r}$ centered at $0$ with radius $1/r$.
Assume also that $p_{i}$ and $q_{0}$ are the only punctures in
$\varphi_{i}^{-1}(\bar{B}_{0}^{r})$ and $\psi
_{0}^{-1}(\bar{B}_{0}^{1/r})$, respectively. In this case we say that
{\it the $i$-th puncture of  $S_{1}$ can be sewn with
the $0$-th puncture of $S_{2}$}. In this case,
 we obtain a sphere with $n+m+1$ punctures and local coordinates by
cutting $\varphi_{i}^{-1}(B_{0}^{r})$ and $\psi_{0}^{-1}(B_{0}^{1/r})$
{}from $S_{1}$ and $S_{2}$, respectively, and then identifying the boundaries
of the resulting surfaces using the map $\varphi_{i} \circ \gamma \circ
\psi_{0}^{-1}$ where $\gamma$ is the map {}from $\Bbb{C}\setminus \{ 0\}$ to
itself defined by $\gamma(z)=1/z$. The punctures (with ordering) of this sphere
with punctures and local coordinates are $p_{0}$, $\dots$, $p_{i-1}$, $q_{1}$,
$\dots$,
$q_{n}$, $p_{i+1}$, $\dots$, $p_{m}$. The local coordinates
vanishing at these punctures are given in the obvious way. Thus we have a
partial
operation. Given two such spheres with punctures and local coordinates, $S_{1}$
and $S_{2}$, with the same number of punctures,
if there is a analytic isomorphism {}from the underlying sphere of $S_{1}$
to the underlying sphere of $S_{2}$ such that the ordered
punctures of $S_{1}$ are
mapped to the ordered punctures of $S_{2}$ and the germs containing
the pull-backs of the
local coordinates of $S_{2}$ are the same as the germs containing
the local coordinates of $S_{1}$, we say that $S_{1}$ and $S_{2}$
are {\it conformally equivalent}. This is an equivalence relation.
The space of conformal equivalence classes of such spheres with punctures and
local coordinates is
called the {\it moduli space of spheres with punctures and local coordinates}.

The moduli space of spheres
with $n+1$ punctures and local coordinates ($n\ge 1$) can be
identified with $K(n)=M^{n-1}\times H \times H_{c}^{n}$ where
$H$ is the set of all sequences $A$ of complex numbers such that
$\mbox{\rm exp}(\displaystyle
\sum_{j=1}^{\infty}A_{j}x^{j+1}\frac{d}{dx})\cdot x$ is a
convergent power series in some neighborhood of $0$,
$H_{c}=\Bbb{C}^{\times} \times H$,
and $M^{n-1}$ is the subset of elements in $\Bbb{C}^{n-1}$ with nonzero
and distinct components. The moduli space of spheres with one punctures and
local coordinates can be
identified
with $K(0)=\{B\in H\;|\;B_{1}=0\}$. Then the moduli space of spheres with
punctures and local coordinates can be identified with
$\cup_{n=1}^{\infty}K(n)$. {}From now on we will refer to $K(n)$,
$n\in \Bbb{N}$ as the moduli
space of spheres with $n+1$ punctures and local coordinates.
The sewing operation for
spheres with punctures and local coordinates induces a partial operation on
$\cup_{n=1}^{\infty}K(n)$.
It is still called the {\it sewing operation} and is denoted
$_{^{i}}\infty_{^{0}}$.
Note that there is an obvious action of $S_{n}$ on $K(n)$ by permuting the
ordering of the  $n$ positively oriented punctures and local coordinates.

Now we have a sequences of sets $K=\{K_{n}\}_{n=1}^{\infty}$
together with partial operations
$_{^{i}}\infty_{^{0}}: K(j)\times K(k)\to K(j+k-1)$,
$j\in \Bbb{Z}_{+}$, $k\in \Bbb{N}$, $i\in \Bbb{Z}_{+}$
and actions of $S_{n}$ on $K(n)$, $n\in \Bbb{Z}_{+}$, respectively. It is
easy to show that the sewing operations satisfy the following conditions
when the sewing operations appearing in the equations below
exist:
\begin{enumerate}

\item For any $j\in \Bbb{Z}_{+}$, $k, l\in \Bbb{N}$, $i_{1}$,
$1\le i_{1}\le j$, $i_{2}$,
$1\le i_{2}\le j+k-1$, $Q_{1}\in K(j)$, $Q_{2}\in K(k)$, $Q_{3}\in K(l)$,
\begin{equation}
(Q_{1\;^{i_{1}}}\infty_{^{0}}Q_{2})_{^{i_{2}}}\infty_{^{0}}Q_{3}=
\left\{\begin{array}{ll}
(Q_{1\;^{i_{2}}}\infty_{^{0}}Q_{3})_{^{l+i_{1}-1}}\infty_{^{0}}Q_{2},
&i_{2}< i_{1},\\
Q_{1\;^{i_{1}}}\infty_{^{0}}(Q_{2\;^{i_{2}-i_{1}+1}}\infty_{^{0}}Q_{3}),
&i_{1}\le i_{2}<i_{1}+k,\\
(Q_{1\;^{i_{2}-j+1}}\infty_{^{0}}Q_{3})_{^{i_{1}}}\infty_{^{0}}Q_{2},
&i_{1}+k\le i_{2}.\end{array} \right.
\end{equation}

\item For any $j\in \Bbb{Z}_{+}$, $k\in \Bbb{N}$, $i$, $1\le i\le k$,
$Q_{1}\in K(j)$, $Q_{2}\in
K(k)$, $\sigma\in S_{j}$ and $\tau\in S_{k}$,
\begin{equation}
\sigma(Q_{1})_{^{i}}\infty_{^{0}}Q_{2}
=\sigma(\stackrel{j-1}{\overbrace{1, \dots, 1}}, k,
\stackrel{j-k}
{\overbrace{1, \dots, 1}})(Q_{1\;^{\sigma(i)}}\infty_{^{0}}Q_{2}),
\end{equation}
\begin{equation}
Q_{1\;^{i}}\infty_{^{0}}\tau(Q_{2})
=(\stackrel{k-1}{\overbrace{1\oplus \cdots \oplus 1}} \oplus
\tau \oplus \stackrel{j-k}{\overbrace{1\oplus \cdots \oplus 1}})
(Q_{1\;^{i}}\infty_{^{0}}Q_{2}).
\end{equation}

\item Let $I=(\bold{0}, (1, \bold{0}))\in H\times (\Bbb{C}^{\times}\times
H)=K(1)$. Then for any $k\in \Bbb{N}$, $i$, $1\le i\le k$, $Q\in K(k)$,
\begin{equation}
Q_{^{i}}\infty_{^{0}}I=I_{^{1}}\infty_{^{0}}Q=Q.
\end{equation}

\end{enumerate}

A sequence $\{\cal{X}(j)\}_{j\in \Bbb{N}}$ of sets equipped with
$\circ_{i}: \cal{X}(j)\times \cal{X}(k)\to \cal{X}(j+k-1)$,
$j\in \Bbb{Z}_{+}$, $k\in \Bbb{N}$, $1\le i\le k$,
actions of $S_{n}$ on $\cal{X}(n)$, $n\in \Bbb{Z}_{+}$, respectively,
and $I\in \cal{X}(1)$ satisfying the conditions (1)--(3) above with
$K(n)$, $n\in \Bbb{N}$, replaced by $\cal{X}(n)$, $_{^{i}}\infty_{^{0}}$
by $\circ_{i}$, is called an {\it operad} \cite{M1}.
If the operations $\circ_{i}$
are only partial and conditions (1)--(3) are satisfied when the operations
in the equations in (1)--(3) exist, it is called a {\it partial operad}
\cite{M2} \cite{HL2} \cite{HL3}.
Thus we see that $K$ is a partial operad. We can also give a topological
structure and a complex analytic structure to $K$ such that
the sewing operations
$_{^{i}}\infty_{^{0}}$ are continuous and complex analytic.

We shall define a (geometric)
vertex
operator algebra to be a ``linear projective representation''
of this partial operad satisfying some additional conditions.
In the representation theory of groups, a linear projective representation
of a group is a linear representation of a central extension of the group.
For $K$, we also have certain extensions which are analogues of
central extensions of groups. These extensions are constructed using
determinant lines  over  spheres with analytically
parametrized boundary.

We describe briefly Segal's work on determinant lines over
Riemann surfaces with analytically
parametrized boundary here.
For details, see \cite{S}. Let $\Sigma$  be a compact Riemann surface
with analytically
parametrized and oriented boundary components. We have the Cauchy-Riemann
operator $\overline{\partial}$ {}from the space $\Omega^{0}(\Sigma)$
of smooth functions on
the surface to the space $\Omega^{0, 1}(\Sigma)$
of $(0, 1)$-forms on the surface. The boundary of $\Sigma$ can be
decomposed as $\partial \Sigma=\cup_{i=1}^{k}C_{i}^{\epsilon_{i}}$ where
for any $i$, $1\le i\le k$,
$C_{i}^{\epsilon_{i}}$ is a connected component of $\partial \Sigma$
and thus is parametrized by an analytic map {}from the circle $S^{1}$ to
$C_{i}^{\epsilon}$ and where $\epsilon_{i}=\pm$ indicates the orientation
of the component. Any smooth function on $C_{i}^{\epsilon_{i}}$ can
be decomposed as the sum of two smooth functions, one of which,
 as a function on
$S^{1}$, has a Fourier expansion of the form $\sum_{n\ge 0}a_{n}e^{2\pi
n\theta i}$ ($\theta$ is the usual parametrization of the circle by angles)
and the other of which, as a function on $S^{1}$, has a
 Fourier expansion of the
form $\sum_{n< 0}a_{n}e^{2\pi
n\theta i}$. If $\epsilon_{i}=+$ ($\epsilon_{i}=-$), that is,
this component is positively (negatively) oriented,
we denote by $\Omega^{0}_{+}(C_{i}^{\epsilon_{i}})$
 the space of all smooth functions on $C_{i}^{\epsilon_{i}}$ which
as functions on $S^{1}$ have Fourier expansions of the form
$\sum_{n\ge 0}a_{n}e^{2\pi
n\theta i}$ ($\sum_{n< 0}a_{n}e^{2\pi
n\theta i}$) and by $\Omega^{0}_{-}(C_{i}^{\epsilon_{i}})$ the space of smooth
functions on $C_{i}^{\epsilon_{i}}$ which
as functions on $S^{1}$ have Fourier expansions of the form
$\sum_{n< 0}a_{n}e^{2\pi
n\theta i}$ ($\sum_{n\ge 0}a_{n}e^{2\pi
n\theta i}$). Thus the space $\Omega^{0}(\partial \Sigma)$ of all smooth
functions on $\partial \Sigma$ can be decomposed as
$\oplus_{i=1}^{k}(\Omega^{0}_{+}(C_{i}^{\epsilon_{i}})\oplus
\Omega^{0}_{-}(C_{i}^{\epsilon_{i}}))$. Following Segal's notation, let
\begin{equation}
\Omega^{0}_{+}(\partial \Sigma)
=\oplus_{i=1}^{k}\Omega^{0}_{-\epsilon_{i}}(C_{i}^{\epsilon_{i}})
\subset \Omega^{0}(\partial \Sigma).
\end{equation}
(Note that in our notation, it is better to denote this space by
$\Omega^{0}_{-}(\partial \Sigma)$. We denote it by
$\Omega^{0}_{+}(\partial \Sigma)$  so that it agrees with Segal's notation.)
Let $\mbox{\rm pr}$ be the composition of the restriction {}from
$\Omega^{0}(\Sigma)$ to $\Omega^{0}(\partial \Sigma)$ and the
projection
{}from
$\Omega^{0}(\partial \Sigma)$ to $\Omega^{0}_{+}(\partial \Sigma)$.
We have an operator
\begin{equation}
\overline{\partial}\oplus \mbox{\rm pr}:
\Omega^{0}(\Sigma)\to \Omega^{0, 1}(\Sigma)
\oplus \Omega^{0}_{+}(\partial \Sigma).
\end{equation}
Using the theory of elliptic boundary problems on manifolds with
boundaries (see, for example, \cite{Ho}), we can show that
$\overline{\partial}\oplus \mbox{\rm pr}$ can be extended to Fredholm
operators {}from suitable Sobolev spaces on $\Sigma$ to direct sums of
closed subspaces of suitable Sobolev spaces on $\Sigma$ and closed
subspaces of suitable Sobolev spaces on $\partial \Sigma$.  In
addition, the kernels of these extensions are equal to the kernel of
$\overline{\partial}\oplus \mbox{\rm pr}$ and the orthogonal
complements of the images of these extensions are in $\Omega^{0,
1}(\Sigma)
\oplus \Omega^{0}_{+}(\partial \Sigma)$. Thus we can regard the kernel and
cokernel of $\overline{\partial}\oplus \mbox{\rm pr}$ as the kernels and
cokernels of its extensions. Since these extensions are Fredholm, the kernel
and cokernel of $\overline{\partial}\oplus \mbox{\rm pr}$ are
finite-dimensional. The determinant line over $\Sigma$ is defined as
\begin{equation}
\mbox{\rm Det}_{\Sigma}=\mbox{\rm Det}\; (\mbox{\rm Ker}\;
(\overline{\partial}\oplus \mbox{\rm pr}))^{*}\otimes
\mbox{\rm Det}\;\mbox{\rm Coker}\;(\overline{\partial}\oplus \mbox{\rm pr})
\end{equation}
where $\mbox{\rm Det}\; (\mbox{\rm Ker}\;
(\overline{\partial}\oplus \mbox{\rm pr}))^{*}$ and
$\mbox{\rm Det}\;
\mbox{\rm Coker}\;(\overline{\partial}\oplus \mbox{\rm pr})$ are the
highest nonzero exterior powers of $(\mbox{\rm Ker}\;
(\overline{\partial}\oplus \mbox{\rm pr}))^{*}$ and
$\mbox{\rm Coker}\;(\overline{\partial}\oplus \mbox{\rm pr})$, respectively.
The main property of determinant lines over
Riemann surfaces with analytically
parametrized and  oriented boundary components is that if we
sew two such Riemann surfaces, $\Sigma_{1}$ and $\Sigma_{2}$,
 by identifying certain boundary components
on $\Sigma_{1}$ to certain boundary components with
opposite orientations on $\Sigma_{2}$ using the given analytic
parametrizations to obtain another such,
denoted by $\Sigma_{1}\infty \Sigma_{2}$,
then there exists a canonical isomorphism
\begin{equation}\label{3.8}
\ell_{\Sigma_{1}, \Sigma_{2}}:
\mbox{\rm Det}_{\Sigma_{1}}\otimes \mbox{\rm Det}_{\Sigma_{2}}
\to \mbox{\rm Det}_{\Sigma_{1}\infty\Sigma_{2}}.
\end{equation}
These determinant lines give a holomorphic line bundle over the moduli space
of Riemann surfaces with oriented and analytically parametrized boundaries,
and there is a canonical connection on this line bundle.
See \cite{S} for more details.

Now we want to use Segal's work described above to define the determinant line
for an element $Q$ of $K$. We need to find a sphere with analytically
parametrized and oriented boundary $\Sigma_{Q}$
determined uniquely by $Q$. For any $Q\in K$,
there is a unique sphere with punctures
and local coordinates in $Q$ such that its underlying sphere is
$\Bbb{C}\cup \{\infty \}$, the negatively oriented puncture is
$\infty$, the last positively oriented puncture is $0$, the value at
$\infty$ of the derivative of the local coordinate map at $\infty$ is
$1$ and all the local coordinate neighborhoods at the punctures
are the preimages under the
local coordinate maps of the maximal open disks (possibly with
infinite radius) centered at $0$ on which the inverses of local
coordinate maps have well-defined analytic extensions.
For any positive real number $r$ and any puncture, consider the closed disk of
radius equal to $r$ times the minimum of $1$ and half of
the radius of the maximal disk above at the puncture.
(To avoid closed disks with
infinite radius, we choose the minimum of $1$ and half of
the radius of the maximal disk instead of half of
the radius of the maximal disk.) For a fixed $r$, a closed disk above is
called a {\it closed disk associated to $r$}. Let $X$ be the set of all
positive real numbers such that if $r\in X$, then at any puncture
 the closed disk associated to $r$ is contained in
the maximal open disk above and preimages under local coordinate maps
of closed
disks associated to $r$ at different punctures
do not intersect each other.
Let $r_{0}=\sup X$ and
$r_{1}=\min(1, \frac{r_{0}}{2})$. (To make sure that $r_{1}$ is not
 $\infty$, we define $r_{1}$ to be $\min(1, \frac{r_{0}}{2})$ instead of
$\frac{r_{0}}{2}$.)
We obtain a Riemann surface with oriented and analytically parametrized
boundary components $\Sigma_{Q}$  by cutting the preimages
of the closed disks associated to $r_{1}$ and giving its
boundary components the obvious orientations and analytic
parametrizations (by first mapping the unit circle to the circle with
radius $r_{1}$).  We define
\begin{equation}
\mbox{\rm Det}_{Q}=\mbox{\rm Det}_{\Sigma_{Q}}.
\end{equation}

We consider annuli on the complex plane with two circles centered at $0$
as boundaries. They can be degenerate in the sense that the two
boundary circles are the
same. With the obvious orientations of the boundary circles and with
multiplications by the radii of the
boundary circles as analytic boundary
parametrizations, these annuli become Riemann surfaces with
oriented and analytically parametrized
boundary components.
If the inner circle of such an annulus is the unit circle, we call it a
{\it canonical annlus with
oriented and analytically parametrized
boundary components}.
For $m, n\in \Bbb{N}$, $Q_{1}\in K(m)$ and $Q_{2}\in K(n)$ such
that $Q_{1^{i}}\infty_{^{0}}Q_{2}$ exists, we can find unique
Riemann surfaces with
oriented and analytically parametrized
boundary components $A$, $B$, $C$, $D$ and $E$
which in general are not connected
and are disjoint unions of
canonical annuli with
oriented and analytically parametrized
boundary components, such that
$((A\infty \Sigma_{Q_{1}})\infty B)\infty(\Sigma_{Q_{2}}\infty C)$
is conformally equivalent to $(D\infty\Sigma_{Q_{1^{i}}\infty_{^{0}}Q_{2}})
\infty E$.
So we have a canonical isomorphism from
$\mbox{\rm Det}_{((A\infty \Sigma_{Q_{1}})\infty B)
\infty(\Sigma_{Q_{2}}\infty C)}$ to
$\mbox{\rm Det}_{(D\infty\Sigma_{Q_{1^{i}}\infty_{^{0}}Q_{2}})
\infty E}$.
It can be shown easily that
$\mbox{\rm Det}_{A}$, $\mbox{\rm Det}_{B}$, $\mbox{\rm Det}_{C}$,
$\mbox{\rm Det}_{D}$ and $\mbox{\rm Det}_{E}$ are
canonically isomorphic to $\Bbb{C}$.
Thus we obtain canonical isomorphisms from
$\mbox{\rm Det}_{\Sigma_{Q_{1}}}\otimes \mbox{\rm Det}_{\Sigma_{Q_{2}}}$
to $\mbox{\rm Det}_{A}\otimes \mbox{\rm Det}_{\Sigma_{Q_{1}}}\otimes
\mbox{\rm Det}_{B}\otimes
\mbox{\rm Det}_{\Sigma_{Q_{2}}}\otimes \mbox{\rm Det}_{C}$ and from
$\mbox{\rm Det}_{D}\otimes \mbox{\rm Det}_{Q_{1^{i}}\infty_{^{0}}Q_{2}}
\otimes \mbox{\rm Det}_{E}$ to $\mbox{\rm Det}_{Q_{1^{i}}\infty_{^{0}}Q_{2}}$.
Composing in the obvious order the three canonical isomorphisms above with
$$\ell_{(A\infty \Sigma_{Q_{1}})\infty B,
\Sigma_{Q_{2}}\infty C}\circ
(1\otimes \ell_{Q_{2}, C})\circ
(\ell_{A\infty \Sigma_{Q_{1}}, B}\otimes 1\otimes 1)\circ
(\ell_{A, \Sigma_{Q_{1}}}\otimes 1\otimes 1\otimes 1),$$
we obtain
a canonical isomorphism
\begin{equation}
\ell^{i}_{Q_{1}, Q_{2}}: \mbox{\rm Det}_{Q_{1}}\otimes \mbox{\rm Det}_{Q_{2}}
\to \mbox{\rm Det}_{Q_{1^{i}}\infty_{^{0}}Q_{2}}.
\end{equation}

Let
\begin{eqnarray}
\tilde{K}(n)&=&\cup_{Q\in  K(n)}\mbox{\rm Det}_{Q}, \;\;\;n\in \Bbb{N},\\
\tilde{K}&=&\{\tilde{K}(n)\}_{n\in \Bbb{N}}.
\end{eqnarray}
Then $\tilde{K}(n)$, $n\in \Bbb{N}$, are holomorphic line bundles (in a
suitable sense) over $K(n)$.
There are also operations in $\tilde{K}$ obtained {}from the sewing operations
in $K$ and the canonical isomorphisms for determinant lines defined as follows:
Let $m, n\in \Bbb{N}$, $i$ an integer satisfying $1\le i\le m$, $Q_{1}\in
K(m)$,
$Q_{2}\in K(n)$,
$\tilde{Q}_{1}\in \mbox{\rm Det}_{Q_{1}} \subset \tilde{K}(m)$ and
$\tilde{Q}_{2}\in \mbox{\rm Det}_{Q_{2}}\subset \tilde{K}(n)$,
such that $Q_{1^{i}}\infty_{^{0}}Q_{2}$ exists.
 We define
\begin{equation}
\tilde{Q}_{1^{i}}\widetilde{\infty}_{^{0}}\tilde{Q}_{2}
=\ell^{i}_{Q_{1}, Q_{2}}(\tilde{Q}_{1} \otimes \tilde{Q}_{2})
\in \mbox{\rm Det}_{Q_{1^{i}}\infty_{^{0}}Q_{2}}\subset \tilde{K}(m+n-1).
\end{equation}
Thus we obtain a partial operation
$_{^{i}}\widetilde{\infty}_{^{0}}: \tilde{K}(m)\times \tilde{K}(n)
\to \tilde{K}(m+n-1)$ for any $m, n\in  \Bbb{N}$ and any
integer $i$ satisfying $1\le i\le m$.
Note that the definition of determinant line over an element
$Q\in K(n)$ for any $n\in \Bbb{N}$ does not use the ordering of the
positively oriented punctures of $Q$. Thus for any $\sigma\in S_{n}$,
$\mbox{\rm Det}_{Q}$ is canonically isomorphic to $\mbox{\rm Det}_{\sigma(Q)}$.
We denote this canonical isomorphism by $\varphi^{\sigma}_{Q}$.
For any $\tilde{Q}\in \mbox{\rm Det}_{Q} \subset \tilde{K}(n)$, we define
\begin{equation}
\sigma(\tilde{Q})=\varphi^{\sigma}_{Q}(\tilde{Q})\in \mbox{\rm Det}_{\sigma(Q)}
\subset \tilde{K}(n).
\end{equation}
We obtain an action of $S_{n}$ on $\tilde{K}(n)$.
Let $\tilde{I}$ be the unique element of $\mbox{\rm Det}_{I}$ satisfying
$\ell^{1}_{I, I}(\tilde{I} \otimes \tilde{I})=\tilde{I}$.
Then the sequence $\tilde{K}$ together with the operations
$$_{^{i}}\widetilde{\infty}_{^{0}}: \tilde{K}(m)\times \tilde{K}(n)\to
\tilde{K}(m+n-1),$$
$m, n\in  \Bbb{N}$,  $1\le i\le m$,
the actions of the symmetric groups and $\tilde{I}$ is a partial operad.
Also the operations $_{^{i}}\widetilde{\infty}_{^{0}}$,
$m, n\in  \Bbb{N}$,  $1\le i\le m$, are all continuous and analytic with
respect to the topological and analytic structures on the holomorphic
line bundles $\tilde{K}(n)$ over $K(n)$,
$n\in \Bbb{N}$.

For any $n\in \Bbb{N}$, there is
a canonical connection on the determinant line bundle $\tilde{K}(n)$
induced from the canonical connection on the determinant line bundle over
the moduli space of Riemann surfaces with oriented and analytic parametrized
boundaries.
Using this connection, we can prove that
the determinant line bundle $\tilde{K}(n)$ is trivial.
Thus for any complex number $c$, a $c$-th power of
determinant line bundle $\tilde{K}(n)$  is well defined. Note
that a $c$-th power of $\tilde{K}(n)$
 is the line bundle whose fibers are the
same as those of $\tilde{K}(n)$ and whose transition
functions are equal to certain branches of
the $c$-th powers of the transition functions
of $\tilde{K}(n)$. The existence of a $c$-power of $\tilde{K}(n)$ means that
we can choose the branches of the  $c$-th powers of the transition functions
of $\tilde{K}(n)$  consistently so that they also give a holomorphic
line bundle, a $c$-th power of $\tilde{K}(n)$.
So we see that because  $\tilde{K}(n)$ is trivial,
there is only one $c$-th power of $\tilde{K}(n)$
and it is in fact canonically isomorphic to
$\tilde{K}(n)$.
We denote the $c$-th power of $\tilde{K}(n)$
by $\tilde{K}^{c}(n)$.
Since, as a line bundle over $K(n)$, $\tilde{K}^{c}(n)$ is
canonically isomorphic to $\tilde{K}(n)$, we shall not distinguish between the
elements of $\tilde{K}(n)$ and the elements of $\tilde{K}^{c}(n)$.
In particular, for any element $\tilde{Q}$ of $\tilde{K}^{c}(n)$
there is $Q\in K(n)$ such that $\tilde{Q}$ is in $\mbox{\rm Det}_{Q}$.
The difference between $\tilde{K}^{c}$ and $\tilde{K}$ is
that the canonical isomorphisms for them are different. We can prove that
we can choose values of $\ell^{i}_{Q_{1}, Q_{2}}$ and $\varphi^{\sigma}_{Q}$
raised to the complex power $c$
(denoted by $(\ell^{i}_{Q_{1}, Q_{2}})^{c}$ and $(\varphi^{\sigma}_{Q})^{c}$,
respectively) consistently
for $m, n\in \Bbb{N}$, $1\le i\le m$, $Q_{1}\in K(m)$, $Q, Q_{2}\in K(n)$
and $\sigma\in S_{n}$,
such that $\tilde{K}^{c}=\{\tilde{K}^{c}(n)\}_{n\in \Bbb{N}}$ becomes
a partial operad;  the operations $_{^{i}}\widetilde{\infty}^{c}_{^{0}}$
are defined in the same
way as those for $_{^{i}}\widetilde{\infty}_{^{0}}$ except that
$\ell^{i}_{Q_{1}, Q_{2}}$ is replaced by $(\ell^{i}_{Q_{1}, Q_{2}})^{c}$,
the actions of the symmetric groups are
defined using $(\varphi^{\sigma}_{Q})^{c}$
and the identity element is $\tilde{I}\in
\tilde{K}^{c}(1)$.
The canonical connection on $\tilde{K}(n)$ gives
a canonical connection on $\tilde{K}^{c}(n)$. Beginning with $\tilde{I}$,
we obtain a section $\psi_{1}$ of $\tilde{K}^{c}(1)$ by parallel transport
(this section is in fact not continuous when $c\ne 0$).
Let $J\in K(0)$ be the conformal equivalence class containing the sphere
$\Bbb{C}\cup \{\infty \}$ with the negatively oriented puncture $\infty$ and
the standard local coordinate $w\to w^{-1}$ vanishing at $\infty$ and let
$\tilde{J}$ be any fixed element of $\mbox{\rm Det}_{J}$. Then beginning with
$\tilde{J}$,
we obtain a section $\psi_{0}$ of $\tilde{K}^{c}(0)$ by parallel transport.
Let $P(1)\in K(2)$ be the conformal equivalence class containing
the sphere $\Bbb{C}\cup \{\infty\}$ with the negatively oriented puncture
$\infty$, the positively oriented punctures $1$ and $0$, the standard local
coordinate $w\to w^{-1}$ vanishing at $\infty$, the standard local coordinate
$w\to w-1$ vanishing at $1$ and the standard local coordinate $w\to w$
vanishing
at $0$. Let $\tilde{P}(1)$ be the unique element of $\mbox{\rm Det}_{P(1)}$
such that $(\ell^{1}_{P(1), J})^{c}(\tilde{P}(1)\otimes \tilde{J})=\tilde{I}$.
Beginning with
$\tilde{P}(1)\in \tilde{K}^{c}(2)$ we obtain a section $\psi_{2}$ of
$\tilde{K}^{c}(2)$ by parallel transport. Since $K$ is generated by
$K(0)$, $K(1)$ and $K(2)$ (which means that any element in $K(n)$ for any
$n\in \Bbb{N}$ can be obtained by sewing elements in
$K(0)$, $K(1)$ and $K(2)$), we obtain sections $\psi_{n}$ of
$\tilde{K}^{c}(n)$, $n\in \Bbb{N}$. It can be shown that $\psi_{n}$,
$n\in \Bbb{N}$, are well-defined. Then we have $\{\psi_{n}\}_{n\in \Bbb{N}}$
which is a section of $\tilde{K}^{c}$.

To define a ``linear representation'' of $\tilde{K}^{c}$, we
first have to construct a partial operad {}from a vector space.  Given a
$\Bbb{Z}$-graded vector space $V=\coprod_{n\in \Bbb{Z}}V_{(n)}$ such that
$\dim V_{(n)} <\infty$, we can construct a partial operad
 in the following way (see \cite{HL2} \cite{HL3}): Let
\begin{eqnarray}
\cal{H}_{V}(n)&=&\hom(V^{n}, \overline{V}),\\
\cal{H}_{V}&=&\{\cal{H}_{V}(n)\}_{n=1}^{\infty}
\end{eqnarray}
where $\overline{V}=\prod_{n\in \Bbb{Z}}V_{(n)}$. Let
$P_{n}$, $n\in \Bbb{Z}$, be the projection {}from $\overline{V}$ to $V_{(n)}$.
For $f\in \cal{H}_{V}(m)$, $g\in \cal{H}_{V}(n)$ and $0\le i \le m$,
if for any $v' \in V'$, $v_{1}, \dots, v_{m+n-1} \in V$ the
series
\begin{equation}
\sum_{n\in \Bbb{Z}}\langle v', f(v_{1}, \dots, v_{i-1}, P_{n}(g(v_{i}, \dots,
v_{i+n-1})), v_{i+n}, \dots, v_{m+n-1})\rangle
\end{equation}
(where $\langle\cdot, \cdot \rangle$ denotes the pairing between
$V'$ and $\overline{V}$) converges,
we say that {\it the contraction} $f_{^{\;i}}\!*_{^{0}}g$
{\it exists} and we define the {\it contraction} $f_{^{\;\;i}}\!*_{^{0}}g\in
\cal{H}_{V}(m+n-1)$
using the values of these series. Note that contractions are partial
operations.
The permutation group $S_{n}$ also acts on $\cal{H}_{V}(n)$ in the
obvious way. We also have the inclusion map $I_{V}\in \cal{H}_{V}(1)=
\hom(V, \overline{V})$. The sequence $\cal{H}_{V}$ together with
the contractions, the actions of the symmetric groups and the inclusion map
$I_{V}$, is a partial operad, called the {\it endomorphism partial operad
of $V$}.

Roughly speaking,  a ``geometric vertex operator algebra'' (or
a ``vertex associative algebra'') is
a $\Bbb{Z}$-graded vector space $V$ equipped with a ``homomorphism''
{}from the partial operad $\tilde{K}^{c}$ to the partial operad
$\cal{H}_{V}$ satisfying some
additional natural axioms. Precisely, we have the following:

\begin{defi}
A {\it geometric vertex operator
algebra of central charge $c$} is
a $\Bbb{Z}$-graded vector space $V$ and a map $\Phi:
\tilde{K}^{c}\longrightarrow
\cal{H}_{V}$ such that $\Phi(\tilde{K}^{c}(n))\subset \cal{H}_{V}(n)$
satisfying:

\begin{enumerate}

\item The positive energy axiom: $V_{(n)}=0$ for $n$ sufficiently small.

\item The grading axiom: Let $Q(a)=(\bold{0}, (a, \bold{0}))\in H\times
(\Bbb{C}^{\times}\times H)=K(1)$ (the conformal equivalence class containing
the sphere $\Bbb{C}\cup\{\infty\}$ with the negatively oriented puncture
$\infty$, the positively oriented puncture $0$, the standard local coordinate
$w\to w^{-1}$ vanishing at $\infty$ and the local coordinate $w\to aw$
vanishing at $0$).
Then for
any $n\in \Bbb{Z}$, $v\in V_{(n)}$, $v'\in V'$,
\begin{equation}
\langle v', \Phi(\psi_{1}(Q(a)))(v)\rangle_{V}=
a^{-n}\langle \cdot, \cdot\rangle_{V_{(n)}}
\end{equation}
where
$\langle \cdot, \cdot\rangle_{n}$ is the pairing between $V'$ and $V_{(n)}$
induced {}from the pairing $\langle \cdot, \cdot \rangle_{V}$ between
$V'$ and $\overline{V}$.

\item The permutation axiom: For any
$n\in \Bbb{N}$, $\sigma\in S_{n}$ and $\tilde{Q}\in \tilde{K}^{c}(n)$,
\begin{equation}
\Phi (\sigma(\tilde{Q}))=\sigma(\Phi (\tilde{Q})).
\end{equation}

\item The analyticity axiom:
For any $n\in \Bbb{N}$, let
\begin{equation}
\nu_{n}=\Phi\circ \psi_{n}: K(n)\to \cal{H}_{V}(n).
\end{equation}
Then for any $v'\in
V'$, $v_{1}, \dots, v_{n}\in  V$, $\langle v',
\nu_{n}(\cdot)(v_{1}\otimes \cdots \otimes
v_{n})\rangle$ as a function of $$(z_{1}, \dots, z_{n-1}; A^{(0)},
(a_{0}^{(1)}, A^{(1)}), \dots, (a_{0}^{(n)}, A^{(n)}))\in
K(n)=M^{n-1}\times H\times (\Bbb{C}^{\times}\times H)^{n}$$
is meromorphic in $z_{1}, \dots, z_{n-1}$
with $z_{i}=0$ and $z_{i}=z_{j}$, $i, j=1,
\dots, n-1$, $i\not =j$,
 as
the only possible poles, and is a Laurent polynomial in $a_{0}^{(1)}, \dots,
a_{0}^{(n)}$ and is a polynomial in the components of $A^{(0)}, \dots,
A^{(n)}$. In addition, for fixed $i, j$,
$1\le i<j\le  n$, and
$v_{i}, v_{j}\in V$ there is an
upper bound, independent of $v_{k}$, $k\ne i, j$,
for the order of the pole $z_{i}-z_{j}$ of the function
$\langle v', \nu_{n}(\cdot)(v_{1}\otimes \cdots\otimes v_{i}\otimes \cdots
\otimes v_{j}\otimes \cdots \otimes v_{n})\rangle$.

\item The sewing axiom: For any $m, n\in \Bbb{N}$,
$\tilde{Q}_{1}\in \mbox{\rm Det}_{Q_{1}}\subset \tilde{K}^{c}(m)$ and
$\tilde{Q}_{2}\in \mbox{\rm Det}_{Q_{1}}\subset \tilde{K}^{c}(n)$ such that
$Q_{1^{i}}\infty_{^{0}}Q_{2}$ exists,
$\Phi(\tilde{Q}_{1})_{^{\;i\!\!}}*_{^{0}}\Phi(\tilde{Q}_{2})$ also exists
and
\begin{equation}
\Phi(\tilde{Q}_{1^{i}}\widetilde{\infty}_{^{0}}\tilde{Q}_{2})
=\Phi(\tilde{Q}_{1})_{^{\;i\!\!}}*_{^{0}}\Phi(\tilde{Q}_{2}).
\end{equation}

\end{enumerate}
\end{defi}

 The definition of
homomorphism {}from one geometric
 vertex operator algebra to another of the same rank is clear.
The following theorem (see \cite{H1}--\cite{H4} \cite{H8})
establishes the equivalence between vertex operator
algebras and geometric vertex operator algebras:

\begin{theo}
The category of geometric vertex operator algebras
of rank $c$ is isomorphic to the category of vertex operator algebras of
rank $c$.
\end{theo}

The map $\nu$ in the definition above can also be
constructed algebraically (see \cite{H1} and \cite{H8}).

\section{Vertex operator algebras and conformal field theories}

The rapidly-evolving theory of vertex operator algebras has been starting
to show its power in the study of many problems related to conformal
field theories. It is expected that in the future this theory will play
a more important role in the study of conformal field theories and related
mathematical problems.

Basically, there are two approaches to conformal field theories. One
is the geometric approach. In physics, many models of conformal field
theories are studied using the path integral method.  Starting {}from
the work of Friedan and Shenker \cite{FS}, physicists have realized
the importance of the moduli space of Riemann surfaces with punctures
in the study of conformal field theories.  The basic mathematical work
in the geometric approach is Segal's definition of conformal field
theory using Riemann surfaces with oriented and analytically
parametrized boundary components \cite{S}.  Motivated by the operator
formalism for the theory of free bosons and free fermions, one
closely related formulation of conformal field theories is given by
Vafa \cite{Va} using Riemann surfaces with punctures and local
coordinates vanishing at these punctures, on a physical level of
rigor. The geometric approach has the advantage that it gives
conceptually satisfactory definitions and it also allows one to derive
many important results using geometric intuition. But the main
difficulty that the geometric approach encountered is that it is very
difficult to construct nontrivial examples satisfying all these
geometric axioms and thus also difficult to discover subtle structures
that a conformal field theory might have. On the other hand, beginning
with the seminal work of Belavin, Polyakov and Zamolodchikov
\cite{BPZ} in physics and the works of Borcherds \cite{B}, Frenkel,
Lepowsky and Meurman \cite{FLM2} in mathematics, another approach, the
algebraic one, provides a practical way for both physicists and
mathematicians to construct concrete examples of conformal field
theories, at least at genus zero and genus one.
There are already many examples of conformal field theories
(in the algebraic formulation) constructed {}from Lie algebras,
lattices, Jordan algebras, $\cal{W}$-algebras (certain associative
algebras similar to the universal enveloping algebra of a Lie
algebra). There are also algebraic methods, for example, methods to
construct orbifold theories and coset models, which give more examples
{}from known ones. But the algebraic approach has the disadvantage
that it mostly constructs and studies only the genus-zero and
genus-one theory.  Also the axioms in the algebraic formulations may
at first seem unfamiliar or complicated  (although they are indeed
completely canonical).  It is therefore necessary and important to
establish rigorously the relationship between the algebraic and
geometric approaches. One of the main ingredients in a conformal field
theory is its ``chiral algebra,'' which is a vertex operator algebra.
The geometric interpretation of vertex operator algebras described in
the preceding section can be viewed as a crucial step of the project of
establishing the equivalence between the two approaches and thus
obtaining examples satisfying the geometric axioms {}from the known
examples satisfying the algebraic axioms. Another step in this
direction is Zhu's work \cite{Z} in which he constructed certain
genus-one correlation functions {}from a vertex operator algebra and
its irreducible modules, assuming that the vertex operator algebra
satisfies certain conditions.

Let me end this exposition with the following picture describing the program of
studying conformal field theories and related mathematical problems
using the representation theory of vertex operator algebras:

\vspace{2em}

\begin{tabular}{c}
Elementary mathematical data (lattices, Lie algebras, \\Jordan algebras,
$\cal{W}$-algebras, etc.)\\
$\Downarrow$\\
Vertex operator algebras, modules, intertwining operators\\
$\Downarrow$\\
Modular functors and conformal field theories (in the sense of
Segal)\\
$\Downarrow$\\
Consequences (Verlinde formulas, modular tensor categories,
knot invariants \\
and three-manifold invariants, monstrous moonshine, etc.)
\end{tabular}

\end{document}